# SIX CONJECTURES WHICH GENERALIZE OR ARE RELATED TO ANDRICA'S CONJECTURE

Florentin Smarandache, Ph D
Associate Professor
Chair of Department of Math & Sciences
University of New Mexico
200 College Road
Gallup, NM 87301, USA
E-mail: smarand@unm.edu

Six conjectures on pairs of consecutive primes are listed below together with examples in each case.

1) The equation $p_{n+1}^x - p_n^x = 1$, (1)

where $p_n$ is the $n^{th}$ prime, has a unique solution in between 0.5 and 1. Checking the first 168 prime numbers (less than 1000), one obtains that:

- The maximum occurs, of course, for $n = 1$, i.e.

$$3^x - 2^x = 1, \text{ when } x = 1.$$

- The minimum occurs for $n = 31$, i.e.

$$127^x - 113^x = 1, \text{ when } x = 0.567148... = a_0 \qquad (2)$$

Thus, Andrica's Conjecture

$$A_n = \sqrt{p_{n+1}} - \sqrt{p_n} < 1$$

is generalized to:

2) $B_n = p_{n+1}^a - p_n^a < 1$, where $a < a_0$. (3)

It is remarkable that the minimum $x$ doesn't occur for $11^x - 7^x = 1$ as in Andrica Conjecture's maximum value, but as in example (2) for $a_0 =$ 0.567148… .

Also, the function $B_n$ in (3) is falling asymptotically as $A_n$ in (2) i.e. in Andrica's Conjecture.

Looking at the prime exponential equations solved with a TI-92 Graphing Calculator (approximately: the bigger the prime number gap is, the smaller solution $x$ for the equation (1); for the same gap between two consecutive primes, the larger the primes, the bigger $x$):

$3^x - 2^x = 1$, has the solution $x = 1.000000$.
$5^x - 3^x = 1$, has the solution $x \approx 0.727160$.
$7^x - 5^x = 1$, has the solution $x \approx 0.763203$.
$11^x - 7^x = 1$, has the solution $x \approx 0.599669$.
$13^x - 11^x = 1$, has the solution $x \approx 0.807162$.
$17^x - 13^x = 1$, has the solution $x \approx 0.647855$.
$19^x - 17^x = 1$, has the solution $x \approx 0.826203$.
$29^x - 23^x = 1$, has the solution $x \approx 0.604284$.

$37^x - 31^x = 1$, has the solution $x \approx 0.624992$.
$97^x - 89^x = 1$, has the solution $x \approx 0.638942$.
$127^x - 113^x = 1$, has the solution $x \approx 0.567148$.
$149^x - 139^x = 1$, has the solution $x \approx 0.629722$.
$191^x - 181^x = 1$, has the solution $x \approx 0.643672$.
$223^x - 211^x = 1$, has the solution $x \approx 0.625357$.
$307^x - 293^x = 1$, has the solution $x \approx 0.620871$.
$331^x - 317^x = 1$, has the solution $x \approx 0.624822$.
$497^x - 467^x = 1$, has the solution $x \approx 0.663219$.
$521^x - 509^x = 1$, has the solution $x \approx 0.666917$.
$541^x - 523^x = 1$, has the solution $x \approx 0.616550$.
$751^x - 743^x = 1$, has the solution $x \approx 0.732707$.
$787^x - 773^x = 1$, has the solution $x \approx 0.664972$.
$853^x - 839^x = 1$, has the solution $x \approx 0.668274$.
$877^x - 863^x = 1$, has the solution $x \approx 0.669397$.
$907^x - 887^x = 1$, has the solution $x \approx 0.627848$.
$967^x - 953^x = 1$, has the solution $x \approx 0.673292$.
$997^x - 991^x = 1$, has the solution $x \approx 0.776959$.

If $x > a_0$, the difference of x-powers of consecutive primes is normally greater than 1. Checking more versions:

$3^{0.99} - 2^{0.99} \approx 0.981037$.
$11^{0.99} - 7^{0.99} \approx 3.874270$.
$11^{0.60} - 7^{0.60} \approx 1.001270$.
$11^{0.59} - 7^{0.59} \approx 0.963334$.
$11^{0.55} - 7^{0.55} \approx 0.822980$.
$11^{0.50} - 7^{0.50} \approx 0.670873$.

$389^{0.99} - 383^{0.99} \approx 5.596550$.

$11^{0.599} - 7^{0.599} \approx 0.997426$.
$17^{0.599} - 13^{0.599} \approx 0.810218$.
$37^{0.599} - 31^{0.599} \approx 0.874526$.
$127^{0.599} - 113^{0.599} \approx 1.230100$.

$997^{0.599} - 991^{0.599} \approx 0.225749$

$127^{0.5} - 113^{0.5} \approx 0.639282$

3) $C_n = p_{n+1}^{1/k} - p_n^{1/k} < 2/k$, where $p_n$ is the n-th prime, and $k \geq 2$ is an integer.

$11^{1/2} - 7^{1/2} \approx 0.670873$.
$11^{1/4} - 7^{1/4} \approx 0.1945837251$.

$11^{1/5} - 7^{1/5} \approx 0.1396211046$ .
$127^{1/5} - 113^{1/5} \approx 0.060837$ .
$3^{1/2} - 2^{1/2} \approx 0.317837$ .
$3^{1/3} - 2^{1/3} \approx 0.1823285204$ .
$5^{1/3} - 3^{1/3} \approx 0.2677263764$ .
$7^{1/3} - 5^{1/3} \approx 0.2029552361$ .
$11^{1/3} - 7^{1/3} \approx 0.3110489078$ .
$13^{1/3} - 11^{1/3} \approx 0.1273545972$ .
$17^{1/3} - 13^{1/3} \approx 0.2199469029$ .
$37^{1/3} - 31^{1/3} \approx 0.1908411993$
$127^{1/3} - 113^{1/3} \approx 0.191938$ .

4) $D_n = p_{n+1}^a - p_n^a < 1/n$ , (4)
where $a < a_0$ and $n$ big enough, $n = n(a)$, holds for infinitely many consecutive primes.

a) Is this still available for $a < 1$?
b) Is there any rank $n_0$ depending on $a$ and $n$ such that (4) is verified for all $n \geq n_0$?

A few examples:

$5^{0.8} - 3^{0.8} \approx 0.21567$ .
$7^{0.8} - 5^{0.8} \approx 1.11938$ .
$11^{0.8} - 7^{0.8} \approx 2.06621$.
$127^{0.8} - 113^{0.8} \approx 4.29973$ .
$307^{0.8} - 293^{0.8} \approx 3.57934$ .
$997^{0.8} - 991^{0.8} \approx 1.20716$ .

5) $p_{n+1} / p_n \leq 5/3$, (5)
the maximum occurs at $n = 2$ .

{The ratio of two consecutive primes is limited,
while the difference $p_{n+1} - p_n$ can be as big as we want!}

6) However, $1/p_n - 1/p_{n+1} \leq 1/6$, and the maximum occurs for n = 1.